\documentclass[reqno,12pt]{amsart}
\usepackage{geometry}
\usepackage{amsmath}
\usepackage{cite}
\usepackage{amsfonts}
\usepackage{mathrsfs}
\usepackage{diagbox}
\usepackage{float}
\usepackage[marginal]{footmisc}

\usepackage{graphicx}
\numberwithin{equation}{section}
\allowdisplaybreaks[4]
\newtheorem{theorem}{Theorem}[section]

\newtheorem{lemma}[theorem]{Lemma}

\begin{document}

\def\square{\hfill${\vcenter{\vbox{\hrule height.4pt \hbox{\vrule width.4pt
height7pt \kern7pt \vrule width.4pt} \hrule height.4pt}}}$}

\title[]
{\small On mean curvature integrals of the outer parallel convex body of constant width}


\author{Zezhen Sun$^{1,\ast}$}

\address{Zezhen Sun$^{1}$, $^{1}$School of Mathematical Sciences, East China Normal University, Shanghai, 200241, China}

\email{52205500017@stu.ecnu.edu.cn}



\thanks{$^{\ast}$Corresponding author. Email address: 52205500017@stu.ecnu.edu.cn}

\begin{abstract}
In this paper, we obtain some results about the mean curvature integrals of the outer parallel convex body of constant width.
\end{abstract}

\maketitle
\section{Introduction}
The mean curvature integrals is a basic concept in Integral Geometry. It connects many geometric invariants, such as area, the degree of the spherical Gauss map, the Euler-Poincar$\acute{e}$ characteristic, the Gauss-Kronecker curvature and so on. Also it is closely related to the Minkowski quermassintegral of convex body and plays an important role in the theory of convex body.

Under the assumptions that $\mathbb{R}^{n}$ is the $n$-dimensional Euclidean space and $L_{r[O]}$ is an $r$-dimensional linear subspace through a fixed point $O$, Santal$\acute{o}$ \cite{san} studied the $l$th mean curvature integral $M^{(n)}_{l}$ of a flattened convex body $K$ in $\mathbb{R}^{n}$ and established the expression of $M^{(n)}_{l}$ in terms of $M^{(r)}_{i}$, where $M^{(r)}_{i}$ is the $i$th mean curvature integral of $K$ in $L_{r[O]}$. Later, Zhou-Jiang \cite{zj}, Jiang-Zeng \cite{jz} and Zeng-Ma-Xia \cite{zmx}
studied the mean curvature integrals $M^{(n)}_{l}$ of the outer parallel convex body $K_{\rho}$, which generalized the results of Santal$\acute{o}$.

Let $\Phi$ be a convex body of constant width $h$ with $C^{2}$ boundary $\partial\Phi$ in $\mathbb{R}^{n}$(any two parallel support hyperplanes of $\Phi$ are always separated by
a constant $h$). Let $(\Phi^{'}_{r})^{(n)}_{\rho}$ be the outer parallel body of $\Phi^{'}_{r}$ in the distance $\rho$ in $\mathbb{R}^{n}$, where $\Phi^{'}_{r}$ is the orthogonal projection of $\Phi$ on the $r$-dimensional linear subspace $L_{r[O]}\subseteq\mathbb{R}^{n}$. Denote by $M^{(n)}_{l}\big(\partial(\Phi^{'}_{r})^{(n)}_{\rho}\big)\,(l=0,1,\cdots,n-1)$ the mean curvature integrals of $(\Phi^{'}_{r})^{(n)}_{\rho}$ and by $M^{(r)}_{l}(\partial\Phi^{'}_{r})\,(l=0,1,\cdots,r-1)$ the  mean curvature integrals of $\Phi^{'}_{r}$ in $L_{r[O]}$. In this paper, we shall prove the following results.

\begin{theorem}\label{them1}
(1)If $l\ge n-r$, then
\begin{align}\label{case1}
 M^{(n)}_{l}\big(\partial(\Phi^{'}_{r})^{(n)}_{\rho}\big)
 &=\sum^{n-l-1}_{j=0}\sum^{n-l-j-1}_{i=0}(-1)^{i}
 \frac{\tbinom{n-l-1}{j}\tbinom{n-l-j-1}{i}\tbinom{r-1}{r-i-1}}{\tbinom{n-1}{n-i-1}}\tfrac{O_{n-i-1}}{O_{r-i-1}}
 M^{(r)}_{r-i-1}(\partial\Phi^{'}_{r})\rho^{j}h^{n-l-j-i-1}.
\end{align}
(2)If $l=n-r-1$, then
\begin{align}\label{case2}
 M^{(n)}_{l}\big(\partial(\Phi^{'}_{r})^{(n)}_{\rho}\big)&=(-1)^{r}\tbinom{n-1}{n-r-1}^{-1}O_{n-r-1}V_{r}(\Phi^{'}_{r})+\sum^{r-1}_{i=0}(-1)^{i}
 \frac{\tbinom{r}{i}\tbinom{r-1}{r-i-1}}{\tbinom{n-1}{n-i-1}}\tfrac{O_{n-i-1}}{O_{r-i-1}}M^{(r)}_{r-i-1}(\partial\Phi^{'}_{r})h^{r-i}\notag\\
 &+\sum^{r}_{j=1}\sum^{r-j}_{i=0}(-1)^{i}\frac{\tbinom{r}{j}\tbinom{r-j}{i}\tbinom{r-1}{r-i-1}}{\tbinom{n-1}{n-i-1}}
 \tfrac{O_{n-i-1}}{O_{r-i-1}}M^{(r)}_{r-i-1}(\partial\Phi^{'}_{r})\rho^{j}h^{r-j-i}.
\end{align}
(3)If $l<n-r-1$, then
\begin{align}\label{case3}
  M^{(n)}_{l}&\big(\partial(\Phi^{'}_{r})^{(n)}_{\rho}\big)\notag\\
  &=\sum^{n-r-l-1}_{j=0}(-1)^{r}
  \tfrac{\tbinom{n-l-1}{j}\tbinom{n-l-j-1}{r}}{\tbinom{n-1}{n-r-1}}O_{n-r-1}V_{r}(\Phi^{'}_{r})\rho^{j}h^{n-l-j-r-1}\notag\\
  &+\sum^{n-r-l}_{j=0}\sum^{r-1}_{i=0}(-1)^{i}
  \frac{\tbinom{n-l-1}{j}\tbinom{n-l-j-1}{i}\tbinom{r-1}{r-i-1}}{\tbinom{n-1}{n-i-1}}\tfrac{O_{n-i-1}}{O_{r-i-1}}
 M^{(r)}_{r-i-1}(\partial\Phi^{'}_{r})\rho^{j}h^{n-l-j-i-1}\notag\\
 &+\sum^{n-l-1}_{j=n-r-l+1}\sum^{n-l-j-1}_{i=0}(-1)^{i}
  \frac{\tbinom{n-l-1}{j}\tbinom{n-l-j-1}{i}\tbinom{r-1}{r-i-1}}{\tbinom{n-1}{n-i-1}}\tfrac{O_{n-i-1}}{O_{r-i-1}}
 M^{(r)}_{r-i-1}(\partial\Phi^{'}_{r})\rho^{j}h^{n-l-j-i-1},
\end{align}
where $V_{r}(\Phi^{'}_{r})$ denotes the $r$-dimensional volume of $\Phi^{'}_{r}$.
\end{theorem}

Based on Theorem \eqref{them1}, we continue to calculate the integral of $M^{(n)}_{l}\big(\partial(\Phi^{'}_{r})^{(n)}_{\rho}\big)$ on Grassmann manifold $G_{r,n-r}$ and get the following result.
\begin{theorem}\label{them2}
(1)If $l\ge n-r$, then
\begin{align*}
&\int_{G_{r,n-r}}M^{(n)}_{l}\big(\partial(\Phi^{'}_{r})^{(n)}_{\rho}\big)dL_{r[O]}
 =\\
 &=\sum^{n-l-1}_{j=0}\sum^{n-l-j-1}_{i=0}(-1)^{i}
 \tfrac{\tbinom{n-l-1}{j}\tbinom{n-l-j-1}{i}\tbinom{r-1}{r-i-1}}{\tbinom{n-1}{n-i-1}}\tfrac{O_{n-i-1}O_{n-2}\cdots O_{n-r}}{O_{r-i-1}O_{r-2}\cdots O_{0}}
 \rho^{j}h^{n-l-j-i-1}M^{(n)}_{n-i-1}(\partial\Phi).
 \end{align*}
(2)If $l=n-r-1$, then
\begin{align}\label{d3}
\int_{G_{r,n-r}}M^{(n)}_{l}&\big(\partial(\Phi^{'}_{r})^{(n)}_{\rho}\big)dL_{r[O]}=\notag\\
&(-1)^{r}\tbinom{n-1}{n-r-1}^{-1}\frac{O_{n-r-1}O_{n-2}\cdots O_{n-r}}{rO_{r-2}\cdots O_{0}}M^{(n)}_{n-r-1}(\partial\Phi)\notag\\
&+\sum^{r-1}_{i=0}(-1)^{i}
 \frac{\tbinom{r}{i}\tbinom{r-1}{r-i-1}}{\tbinom{n-1}{n-i-1}}\tfrac{O_{n-i-1}O_{n-2}\cdots O_{n-r}}{O_{r-i-1}O_{r-2}\cdots O_{0}}h^{r-i}M^{(n)}_{n-i-1}(\partial\Phi)\notag\\
 &+\sum^{r}_{j=1}\sum^{r-j}_{i=0}(-1)^{i}\frac{\tbinom{r}{j}\tbinom{r-j}{i}\tbinom{r-1}{r-i-1}}{\tbinom{n-1}{n-i-1}}
 \tfrac{O_{n-i-1}O_{n-2}\cdots O_{n-r}}{O_{r-i-1}O_{r-2}\cdots O_{0}}\rho^{j}h^{r-j-i}M^{(n)}_{n-i-1}(\partial\Phi).
\end{align}
(3)If $l<n-r-1$, then
\begin{align*}
&\int_{G_{r,n-r}}M^{(n)}_{l}\big(\partial(\Phi^{'}_{r})^{(n)}_{\rho}\big)dL_{r[O]}\notag=\\
  &=\sum^{n-r-l-1}_{j=0}(-1)^{r}
 \tfrac{\tbinom{n-l-1}{j}\tbinom{n-l-j-1}{r}}{\tbinom{n-1}{n-r-1}}\tfrac{O_{n-r-1}O_{n-2}\cdots O_{n-r}}{rO_{r-2}\cdots O_{0}}\rho^{j}h^{n-l-j-r-1}M^{(n)}_{n-r-1}(\partial\Phi)\notag\\
  &+\sum^{n-r-l}_{j=0}\sum^{r-1}_{i=0}(-1)^{i}
  \frac{\tbinom{n-l-1}{j}\tbinom{n-l-j-1}{i}\tbinom{r-1}{r-i-1}}{\tbinom{n-1}{n-i-1}}\tfrac{O_{n-i-1}O_{n-2}\cdots O_{n-r}}{O_{r-i-1}O_{r-2}\cdots O_{0}}\rho^{j}h^{n-l-j-i-1}
 M^{(n)}_{n-i-1}(\partial\Phi)\notag\\
 &+\sum^{n-l-1}_{j=n-r-l+1}\sum^{n-l-j-1}_{i=0}(-1)^{i}
  \frac{\tbinom{n-l-1}{j}\tbinom{n-l-j-1}{i}\tbinom{r-1}{r-i-1}}{\tbinom{n-1}{n-i-1}}\tfrac{O_{n-i-1}O_{n-2}\cdots O_{n-r}}{O_{r-i-1}O_{r-2}\cdots O_{0}}\rho^{j}h^{n-l-j-i-1}
 M^{(n)}_{n-i-1}(\partial\Phi).
\end{align*}
\end{theorem}

\section{Preliminaries}
A set in the Euclidean space $\mathbb{R}^{n}$ is called convex if and only if it contains, with each pair
of its points, the entire line segment joining them. A convex set with nonempty interior is
called a convex body. The boundary $\partial K$ of a convex body $K$ is called a convex hypersurface.

Let $K$ be a convex body and $O$ be a fixed point in $\mathbb{R}^{n}$. Consider all the $(n-r)$-dimensional planes $L_{n-r[O]}$ through $O$.  The orthogonal projection $K^{'}_{n-r}$ of $K$ onto $L_{n-r[O]}$ is defined by the convex set of $K$.  Denoted by $V(K^{'}_{n-r})$  the volume of $K^{'}_{n-r}$ and $O_{m}$ the
area of the $m$-dimensional unit sphere. Then the quermassintegral all the intersection points of  $W^{(n)}_{r}(K)$ of $K$ is defined by (see \cite{san2})
\begin{equation}\label{e0}
 W^{(n)}_{r}(K)=\frac{(n-r)O_{n-1}}{nO_{n-r-1}}E\big(V(K^{'}_{n-r})\big)=\frac{(n-r)O_{r-1}\cdots O_{0}}{nO_{n-2}\cdots O_{n-r-1}}I_{r}(K),\,\,r=1,\cdots,n-1,
\end{equation}
where
\begin{equation}\label{e1}
E\big(V(K^{'}_{n-r})\big)=\frac{I_{r}(K)}{m(G_{n-r,r})},
\end{equation}
\begin{equation}\label{e2}
I_{r}(K)=\int_{G_{n-r,r}}V(K^{'}_{n-r})dL_{n-r[O]}=\int_{G_{r,n-r}}V(K^{'}_{n-r})dL_{r[O]},
\end{equation}
\begin{equation}\label{e3}
 m(G_{n-r,r})=\int_{G_{n-r,r}}dL_{n-r[O]}=\int_{G_{r,n-r}}dL_{r[O]}=\frac{O_{n-1}\cdots O_{n-r}}{O_{r-1}\cdots O_{1}O_{0}}.
\end{equation}
For completeness, we define $W^{(n)}_{0}(K)=I_{0}(K)=V(K),\,W^{(n)}_{n}(K)=O_{n-1}/n$.

Quermassintegral is introduced by Minkowski. It describes the mean value of the projected volumes of a convex body and is a strong tool in the theory of convex body. Cauchy, Kubota, Steiner and so on obtained many famous formulas(see \cite{san2},\cite{ren},\cite{sch}).

The outer parallel body $K_{\rho}$ in the distance $\rho$ of a convex set $K$ is the union of all solid
spheres of radius $\rho$ the centers of which are points of $K$. The boundary $\partial K_{\rho}$ is called the parallel hypersurface of $\partial K$ in the distance $\rho$.  we have the following Steiner formula for the outer parallel body $K_{\rho}$.
\begin{equation}\label{vk}
 V(K_{\rho})=\sum^{n}_{i=0}\begin{pmatrix}
  n\\i
  \end{pmatrix}W^{(n)}_{i}(K)\rho^{i}.
 \end{equation}
As a consequence of the Steiner formula we have
\begin{equation}\label{wnk}
 W^{(n)}_{i}(K_{\rho})=\sum^{n-i}_{j=0}\begin{pmatrix}
 n-i\\j
 \end{pmatrix}W^{(n)}_{i+j}(K)\rho^{j},\,i=0,1,\cdots,n.
\end{equation}
Moreover, we have the relation between the mean curvature integrals of $\partial K$ and the
Minkowski quermassintegrals of $K$
 \begin{equation}\label{wni}
 M^{(n)}_{i}(\partial K)=nW^{(n)}_{i+1}(K),\,i=0,1,\cdots,n-1.
 \end{equation}
Note that the Minkowski quermassintegrals $W^{(n)}_{i}$ are well defined for any convex figure, whereas $W^{(n)}_{i}(\partial K)$ makes sense only if $\partial K$ is of class $C^{2}$.

Let $K$ be a convex body in $\mathbb{R}^{n}$, then $\partial K$ is an $(n-1)$ -dimensional convex hypersurface. Assuming that $\partial K$ is of class $C^{2}$ and $P$ is a point of $\partial K$,  we choose$
e_{1},\cdots,e_{n-1}$ to be the principal curvature directions at the point $P$. Further, we suppose that $\kappa_{1},\cdots,\kappa_{n-1}$ are the principal curvatures at the point $P$, which correspond to the principal curvature directions.
Consider the Gauss map $G:p \rightarrow N(p)$, , whose differential
\begin{equation*}
dG_{p}:x^{'}(t) \rightarrow N^{'}(t)\,\,(x(0)=p)
\end{equation*}
satisfies Rodrigues' equations,
\begin{equation*}
dG_{p}(e_{i})=-\kappa_{i}e_{i},\,\,i=1,\cdots,n-1.
\end{equation*}
Then we have the mean curvature
\begin{equation*}
H=\frac{1}{n-1}(\kappa_{1}+\cdots+\kappa_{n-1})=-\frac{1}{n-1}trace(dG_{p}),
\end{equation*}
and the Gauss-Kronecker curvature
\begin{equation*}
K=\kappa_{1}\cdots\kappa_{n-1}=(-1)^{n-1}det(dG_{p}).
\end{equation*}
The $i$th order mean curvature is the $i$th order elementary symmetric function of the principal curvatures. We denote by $H_{i}$ the $i$th order mean curvature normalized such that
\begin{equation*}
\prod^{n-1}_{i=1}(1+t\kappa_{i})=\sum^{n-1}_{i=0}H_{i}t^{i}.
\end{equation*}
Thus, $H_{1}=H$ is the mean curvature and $H_{n-1}$ is the Gauss-Kronecker curvature.

The $i$th order mean curvature integral $M^{(n)}_{i}$ of $\partial K$ at $P$ is defined by
\begin{equation*}
M^{(n)}_{i}(\partial K)=\int_{\partial K}H_{i}d\sigma=
\begin{pmatrix}
 n-i\\j
 \end{pmatrix}^{-1}\int_{\partial K}{\kappa_{j_{1}},\cdots,\kappa_{j_{i}}}d\sigma,\,\,i=1,\cdots,n-1,
\end{equation*}
where ${\kappa_{j_{1}},\cdots,\kappa_{j_{i}}}$ denotes the $i$th elementary symmetric function of the principal curvatures and $d\sigma$ is the area element of $\partial K$. Let $M^{(n)}_{0}(\partial K)=F$,  the area of $\partial K$, for completeness.

Let $K$ be a convex body in $\mathbb{R}^{n}$ and $L_{r[O]}(r<n)$ $r$-dimensional plane through fixed point $O$ in $\mathbb{R}^{n}$. Denote by $K^{'}_{r}$ the orthogonal projection of $K$ onto $L_{r[O]}$. Now, let
$M^{(r)}_{q}(\partial K^{'}_{r})(q=0,1,\cdots,r-1)$  be the mean curvature integrals of $\partial K^{'}_{r}$
as a convex surface in $L_{r[O]}$ and $M^{(n)}_{q}(\partial K^{'}_{r})(q=0,1,\cdots,r-1)$ the mean curvature integrals of $K^{'}_{r}$ as a flattened convex body in $\mathbb{R}^{n}$, then we have the relations between them obtained by Santal$\acute{o}$ (see \cite{san2},\cite{ren}).
\begin{lemma}
(1)If $q\ge n-r$, then
\begin{equation}\label{sa1}
M^{(n)}_{q}(\partial K^{'}_{r})=\frac{\begin{pmatrix}
 r-1\\q-n+r
 \end{pmatrix}}{\begin{pmatrix}
 n-1\\q
 \end{pmatrix}}\frac{O_{q}}{O_{q-n+r}}M^{(r)}_{q-n+r}(\partial K^{'}_{r}).
\end{equation}
(2)If $q=n-r-1$ ,then
\begin{equation}\label{sa2}
M^{(n)}_{n-r-1}(\partial K^{'}_{r})=\begin{pmatrix}
 n-1\\n-r-1
 \end{pmatrix}^{-1}O_{n-r-1}V_{r}(K^{'}_{r}),
\end{equation}
where $V_{r}(K^{'}_{r})$ denotes the r-dimensional volume of $K^{'}_{r}$.
\\(3)If $q<n-r-1$ ,then
\begin{equation}\label{sa3}
M^{(n)}_{q}(\partial K^{'}_{r})=0.
\end{equation}
\end{lemma}
We also need the following result(see \cite{mar}).
\begin{lemma}
Let $\Phi$ be an convex body of constant width $h$ in $\mathbb{R}^{n}$. Then
\begin{equation}\label{ph}
W^{(n)}_{s}(\Phi)=\sum^{n-s}_{i=0}(-1)^{i}\begin{pmatrix}
 n-s\\i
 \end{pmatrix}W^{(n)}_{n-i}(\Phi)h^{n-s-i},\,\,s=0,1,\cdots,n,
\end{equation}
where $W^{(n)}_{s}(\Phi)$ is the quermassintegral of $\Phi$.
\end{lemma}

\section{The proof of main Theorems}
\textbf{\emph{The proof of  Theorem $\ref{them1}$}}
\begin{proof}
By \eqref{wnk} and \eqref{ph} we have
\begin{align}\label{w0}
 W^{(n)}_{l}(\Phi_{\rho})&=\sum^{n-l}_{j=0}
 \begin{pmatrix}
 n-l\\j
 \end{pmatrix}W^{(n)}_{l+j}(K)\rho^{j}\notag\\
 &=\sum^{n-l}_{j=0}\sum^{n-l-j}_{i=0}(-1)^{i}
 \begin{pmatrix}
 n-l\\j
 \end{pmatrix}
 \begin{pmatrix}
 n-l-j\\i
 \end{pmatrix}W^{(n)}_{n-i}(\Phi)\rho^{j}h^{n-l-j-i},\,\,l=0,1,\cdots,n.
\end{align}
Applying this formula to the convex body $(\Phi^{'}_{r})^{(n)}_{\rho}$, we have
\begin{equation}\label{wa}
 W^{(n)}_{l}\big((\Phi^{'}_{r})^{(n)}_{\rho}\big)
 =\sum^{n-l}_{j=0}\sum^{n-l-j}_{i=0}(-1)^{i}
 \begin{pmatrix}
 n-l\\j
 \end{pmatrix}
 \begin{pmatrix}
 n-l-j\\i
 \end{pmatrix}W^{(n)}_{n-i}(\Phi^{'}_{r})\rho^{j}h^{n-l-j-i},\,\,l=0,1,\cdots,n.
\end{equation}
By \eqref{wni}, we obtain
\begin{equation}\label{wb}
 M^{(n)}_{l}\big(\partial (\Phi^{'}_{r})^{(n)}_{\rho}\big)=nW^{(n)}_{l+1}\big((\Phi^{'}_{r})^{(n)}_{\rho}\big),\,\,l=0,1,\cdots,n-1.
 \end{equation}
We apply \eqref{wa} and \eqref{wb} to get for $l=0,1,\cdots,n$,
\begin{align}\label{wc}
 M^{(n)}_{l}\big(\partial(\Phi^{'}_{r})^{(n)}_{\rho}\big)
 &=n\sum^{n-l-1}_{j=0}\sum^{n-l-j-1}_{i=0}(-1)^{i}
 \begin{pmatrix}
 n-l-1\\j
 \end{pmatrix}
 \begin{pmatrix}
 n-l-j-1\\i
 \end{pmatrix}W^{(n)}_{n-i}(\Phi^{'}_{r})\rho^{j}h^{n-l-j-i-1}\notag\\
 &=\sum^{n-l-1}_{j=0}\sum^{n-l-j-1}_{i=0}(-1)^{i}
 \begin{pmatrix}
 n-l-1\\j
 \end{pmatrix}
 \begin{pmatrix}
 n-l-j-1\\i
 \end{pmatrix}M^{(n)}_{n-i-1}(\partial\Phi^{'}_{r})\rho^{j}h^{n-l-j-i-1}.
\end{align}

Now, we are ready to compute the mean curvature integral of $\partial(\Phi^{'}_{r})^{(n)}_{\rho}$ from the below three cases.

\textbf{Case 1.}If $l\ge n-r$, we have $n-i-1\ge n-(n-l-j-1)-1\ge l\ge n-r$. Then by \eqref{wc} and Santal$\acute{o}$'s result \eqref{sa1}, we can get
\begin{align*}
 M^{(n)}_{l}\big(\partial(\Phi^{'}_{r})^{(n)}_{\rho}\big)
 &=\sum^{n-l-1}_{j=0}\sum^{n-l-j-1}_{i=0}(-1)^{i}
 \tbinom{n-l-1}{j}
 \tbinom{n-l-j-1}{i}
 \frac{\tbinom{r-1}{r-i-1}}{\tbinom{n-1}{n-i-1}}\tfrac{O_{n-i-1}}{O_{r-i-1}}
 M^{(r)}_{r-i-1}(\partial\Phi^{'}_{r})\rho^{j}h^{n-l-j-i-1}.
 \end{align*}

\textbf{Case 2.}If $l=n-r-1$, then \eqref{wc} can be expressed as
\begin{align}
 M^{(n)}_{l}\big(\partial(\Phi^{'}_{r})^{(n)}_{\rho}\big)
 &=\sum^{r}_{j=0}\sum^{r-j}_{i=0}(-1)^{i}
 \tbinom{r}{j}
 \tbinom{r-j}{i}M^{(n)}_{n-i-1}(\partial\Phi^{'}_{r})\rho^{j}h^{r-j-i}.
\end{align}
Note that $n-i-1\ge n-r-1$, thus the right side of the above equation can be decomposed into three parts
\begin{align*}
 M^{(n)}_{l}\big(\partial(\Phi^{'}_{r})^{(n)}_{\rho}\big)
 &=(-1)^{r}M^{(n)}_{n-r-1}(\partial\Phi^{'}_{r})+\sum^{r-1}_{i=0}(-1)^{i}
 \tbinom{r}{i}M^{(n)}_{n-i-1}(\partial\Phi^{'}_{r})h^{r-i}\\
 &+\sum^{r}_{j=1}\sum^{r-j}_{i=0}(-1)^{i}\tbinom{r}{j}
 \tbinom{r-j}{i}M^{(n)}_{n-i-1}(\partial\Phi^{'}_{r})\rho^{j}h^{r-j-i}\\
 &=(-1)^{r}\tbinom{n-1}{n-r-1}^{-1}O_{n-r-1}V_{r}(\Phi^{'}_{r})+\sum^{r-1}_{i=0}(-1)^{i}
 \frac{\tbinom{r}{i}\tbinom{r-1}{r-i-1}}{\tbinom{n-1}{n-i-1}}\tfrac{O_{n-i-1}}{O_{r-i-1}}M^{(r)}_{r-i-1}(\partial\Phi^{'}_{r})h^{r-i}\\
 &+\sum^{r}_{j=1}\sum^{r-j}_{i=0}(-1)^{i}\frac{\tbinom{r}{j}\tbinom{r-j}{i}\tbinom{r-1}{r-i-1}}{\tbinom{n-1}{n-i-1}}
 \tfrac{O_{n-i-1}}{O_{r-i-1}}M^{(r)}_{r-i-1}(\partial\Phi^{'}_{r})\rho^{j}h^{r-j-i}.
\end{align*}
where the last equation follows from \eqref{sa1} and \eqref{sa2}.

\textbf{Case 3.}If $l<n-r-1$, we first analyze the lower index of mean curvature integral $M^{(n)}_{n-i-1}$ in the equation \eqref{wc} and get the following table without considering other coefficients for the time being.

\begin{table}[H]
  \centering
  \fontsize{8}{10}\selectfont
  \caption{$M^{(n)}_{n-i-1}$}
  \begin{tabular}{|c|c|c|c|c|c|c|c|c|c|c|}
  \hline
   \diagbox{$j$}{$i$}&0 &1 & $\cdots\cdots$ &$r-2$ &$r-1$ &r &$r+1$ & $\cdots\cdots$ &$n-l-2$ &$n-l-1$\\
   \hline
    0 &$M^{(n)}_{n-1}$ &$M^{(n)}_{n-2}$ &$\cdots\cdots$ &$M^{(n)}_{n-r+1}$ &$M^{(n)}_{n-r}$ &$M^{(n)}_{n-r-1}$ & $M^{(n)}_{n-r-2}$ & $\cdots\cdots$ & $M^{(n)}_{l+1}$  & $M^{(n)}_{l}$ \\
    \hline
    1 &$M^{(n)}_{n-1}$ &$M^{(n)}_{n-2}$ &$\cdots\cdots$ &$M^{(n)}_{n-r+1}$ &$M^{(n)}_{n-r}$ &$M^{(n)}_{n-r-1}$ & $M^{(n)}_{n-r-2}$ & $\cdots\cdots$ & $M^{(n)}_{l+1}$  & \\
    \hline
    $\cdots\cdots$ &$M^{(n)}_{n-1}$ &$M^{(n)}_{n-2}$ &$\cdots\cdots$ &$M^{(n)}_{n-r+1}$ &$M^{(n)}_{n-r}$ &$\cdots\cdots$ &$\cdots\cdots$ & $\cdots\cdots$ & &\\
    \hline
    $n-r-l-1$ &$M^{(n)}_{n-1}$ & $M^{(n)}_{n-2}$   & $\cdots\cdots$ &$M^{(n)}_{n-r+1}$ &$M^{(n)}_{n-r}$ &$M^{(n)}_{n-r-1}$ & & & &\\
    \hline
    $n-r-l$    &$M^{(n)}_{n-1}$ & $M^{(n)}_{n-2}$   & $\cdots\cdots$ &$M^{(n)}_{n-r+1}$ &$M^{(n)}_{n-r}$ & & & & &\\
    \hline
    $n-r-l+1$  &$M^{(n)}_{n-1}$  & $M^{(n)}_{n-2}$ &$\cdots\cdots$ &$M^{(n)}_{n-r+1}$ & & & & & &\\
    \hline
    $\cdots\cdots$   & $M^{(n)}_{n-1}$ &$M^{(n)}_{n-2}$ &$\cdots\cdots$ & & & & & & &\\
    \hline
   $n-l-2$ &$M^{(n)}_{n-1}$ &$M^{(n)}_{n-2}$ & & & & & & &  &\\
    \hline
    $n-l-1$ &$M^{(n)}_{n-1}$ & & & & & & & & &\\
    \hline
   \end{tabular}
\end{table}
By \eqref{sa3}, we know that when $i>r$, $M^{(n)}_{n-i-1}(\partial\Phi^{'}_{r})=0$. Therefore, combined with the above table,  \eqref{wc} can be rewritten as
\begin{align*}
  M^{(n)}_{l}&\big(\partial(\Phi^{'}_{r})^{(n)}_{\rho}\big)\\
  &=\sum^{n-r-l}_{j=0}\sum^{r-1}_{i=0}(-1)^{i}
  \tbinom{n-l-1}{j}
  \tbinom{n-l-j-1}{i}
  M^{(n)}_{n-i-1}(\partial\Phi^{'}_{r})\rho^{j}h^{n-l-j-i-1}\\
 &+\sum^{n-l-1}_{j=n-r-l+1}\sum^{n-l-j-1}_{i=0}(-1)^{i}
  \tbinom{n-l-1}{j}
  \tbinom{n-l-j-1}{i}
  M^{(n)}_{n-i-1}(\partial\Phi^{'}_{r})\rho^{j}h^{n-l-j-i-1}\\
  &+\sum^{n-r-l-1}_{j=0}(-1)^{r}
  \tbinom{n-l-1}{j}
  \tbinom{n-l-j-1}{r}
  M^{(n)}_{n-r-1}(\partial\Phi^{'}_{r})\rho^{j}h^{n-l-j-r-1}\\
  &=\sum^{n-r-l}_{j=0}\sum^{r-1}_{i=0}(-1)^{i}
  \tbinom{n-l-1}{j}
  \tbinom{n-l-j-1}{i}
  \frac{\tbinom{r-1}{r-i-1}}{\tbinom{n-1}{n-i-1}}\tfrac{O_{n-i-1}}{O_{r-i-1}}
 M^{(r)}_{r-i-1}(\partial\Phi^{'}_{r})\rho^{j}h^{n-l-j-i-1}\\
 &+\sum^{n-l-1}_{j=n-r-l+1}\sum^{n-l-j-1}_{i=0}(-1)^{i}
  \tbinom{n-l-1}{j}
  \tbinom{n-l-j-1}{i}
  \frac{\tbinom{r-1}{r-i-1}}{\tbinom{n-1}{n-i-1}}\tfrac{O_{n-i-1}}{O_{r-i-1}}
 M^{(r)}_{r-i-1}(\partial\Phi^{'}_{r})\rho^{j}h^{n-l-j-i-1}\\
 &+\sum^{n-r-l-1}_{j=0}(-1)^{r}
  \tbinom{n-l-1}{j}
  \tbinom{n-l-j-1}{r}
  \tbinom{n-1}{n-r-1}^{-1}O_{n-r-1}V_{r}(\Phi^{'}_{r})\rho^{j}h^{n-l-j-r-1}
\end{align*}
where the last equation follows from \eqref{sa1} and \eqref{sa2}. This completes the proof Theorem \eqref{them1}.
\end{proof}

\textbf{\emph{The proof of  Theorem $\ref{them2}$}}
\begin{proof}
By \eqref{e1}, we have
\begin{equation*}
I_{n-r}(\Phi)=\int_{G_{r,n-r}}V(\Phi^{'}_{r})dL_{r[O]}=\int_{G_{n-r,r}}V(\Phi^{'}_{r})dL_{n-r[O]}.
\end{equation*}
Divided by $m(G_{r,n-r})$,the volume of Grassmann manifold $G_{r,n-r}$, we obtain the mean value
of the projection volumes $V(\phi^{'}_{r})$
\begin{equation*}
E\big(V(\Phi^{'}_{r})\big)=\frac{I_{n-r}(\Phi)}{m(G_{r,n-r})}=\frac{O_{r-1}\cdots O_{1}O_{0}}{O_{n-1}\cdots O_{n-r}}I_{n-r}(\Phi).
\end{equation*}
Recalling the definition of quermassintegral \eqref{e0}, we have
\begin{equation*}
 W^{(n)}_{n-r}(\Phi)=\frac{rO_{n-1}}{nO_{r-1}}E\big(V(\Phi^{'}_{r})\big)=\frac{rO_{r-2}\cdots O_{0}}{nO_{n-2}\cdots O_{n-r}}I_{n-r}(\Phi),\,\,r=1,\cdots,n-1.
\end{equation*}
Apply the above equation to $\Phi_{\rho}$ and use the fact $(\Phi_{\rho})^{'}_{r}=(\Phi^{'}_{r})^{(r)}_{\rho}$(\cite{ren}) to get
\begin{equation}\label{c1}
 W^{(n)}_{n-r}(\Phi_{\rho})=\frac{rO_{r-2}\cdots O_{0}}{nO_{n-2}\cdots O_{n-r}}\int_{G_{r,n-r}}V\big((\Phi^{'}_{r})^{(r)}_{\rho}\big)dL_{r[O]},\,\,r=1,\cdots,n-1.
\end{equation}
On the other hand, by \eqref{w0}, we get
\begin{equation}\label{c2}
 W^{(n)}_{n-r}(\Phi_{\rho})=\sum^{r}_{j=0}
 \begin{pmatrix}
 r\\j
 \end{pmatrix}W^{(n)}_{n-r+j}(\Phi)\rho^{j},\,\,r=0,1,\cdots,n-1.
\end{equation}
And by Steiner formula \eqref{vk} we have
\begin{equation}\label{c3}
 V\big((\Phi^{'}_{r})^{(r)}_{\rho}\big)=\sum^{r}_{j=0}
 \begin{pmatrix}
 r\\j
 \end{pmatrix}W^{(r)}_{j}(\Phi^{'}_{r})\rho^{j},\,\,r=0,1,\cdots,n-1.
 \end{equation}
Combined with \eqref{c1}, \eqref{c1} and \eqref{c2}, we get
\begin{equation*}
\int_{G_{r,n-r}}W^{(r)}_{j}(\Phi^{'}_{r})dL_{r[O]}=\frac{nO_{n-2}\cdots O_{n-r}}{rO_{r-2}\cdots O_{0}}
W^{(n)}_{n-r+j}(\Phi).
\end{equation*}
Using \eqref{wni} yields
\begin{equation}\label{c4}
\int_{G_{r,n-r}}M^{(r)}_{t}(\partial\Phi^{'}_{r})dL_{r[O]}=\frac{O_{n-2}\cdots O_{n-r}}{O_{r-2}\cdots O_{0}}
M^{(n)}_{n-r+t}(\partial\Phi),\,\,t=0,1,\cdots,r-1.
\end{equation}

Now, we are ready to compute the integral of $M^{(n)}_{l}\big(\partial(\Phi^{'}_{r})^{(n)}_{\rho}\big)$ on Grassmann manifold $G_{r,n-r}$.

\textbf{Case 1.}When $q\ge n-r$, by \eqref{case1} and \eqref{c4}
\begin{align*}
&\int_{G_{r,n-r}}M^{(n)}_{l}\big(\partial(\Phi^{'}_{r})^{(n)}_{\rho}\big)dL_{r[O]}
 =\\
 &\int_{G_{r,n-r}}\sum^{n-l-1}_{j=0}\sum^{n-l-j-1}_{i=0}(-1)^{i}
 \frac{\tbinom{n-l-1}{j}\tbinom{n-l-j-1}{i}\tbinom{r-1}{r-i-1}}{\tbinom{n-1}{n-i-1}}\tfrac{O_{n-i-1}}{O_{r-i-1}}
 M^{(r)}_{r-i-1}(\partial\Phi^{'}_{r})\rho^{j}h^{n-l-j-i-1}dL_{r[O]}\\
 &=\sum^{n-l-1}_{j=0}\sum^{n-l-j-1}_{i=0}(-1)^{i}
 \tfrac{\tbinom{n-l-1}{j}\tbinom{n-l-j-1}{i}\tbinom{r-1}{r-i-1}}{\tbinom{n-1}{n-i-1}}\tfrac{O_{n-i-1}O_{n-2}\cdots O_{n-r}}{O_{r-i-1}O_{r-2}\cdots O_{0}}
 \rho^{j}h^{n-l-j-i-1}M^{(n)}_{n-i-1}(\partial\Phi).
 \end{align*}

\textbf{Case 2.}When $q=n-r-1$, by \eqref{case2} and \eqref{c4}
\begin{align}\label{d1}
\int_{G_{r,n-r}}M^{(n)}_{l}&\big(\partial(\Phi^{'}_{r})^{(n)}_{\rho}\big)dL_{r[O]}=\notag\\
&(-1)^{r}\tbinom{n-1}{n-r-1}^{-1}O_{n-r-1}\int_{G_{r,n-r}}V_{r}(\Phi^{'}_{r})dL_{r[O]}\notag\\
&+\sum^{r-1}_{i=0}(-1)^{i}
 \frac{\tbinom{r}{i}\tbinom{r-1}{r-i-1}}{\tbinom{n-1}{n-i-1}}\tfrac{O_{n-i-1}}{O_{r-i-1}}h^{r-i}\int_{G_{r,n-r}}M^{(r)}_{r-i-1}(\partial\Phi^{'}_{r})dL_{r[O]}\notag\\
 &+\sum^{r}_{j=1}\sum^{r-j}_{i=0}(-1)^{i}\frac{\tbinom{r}{j}\tbinom{r-j}{i}\tbinom{r-1}{r-i-1}}{\tbinom{n-1}{n-i-1}}
 \tfrac{O_{n-i-1}}{O_{r-i-1}}\rho^{j}h^{r-j-i}\int_{G_{r,n-r}}M^{(r)}_{r-i-1}(\partial\Phi^{'}_{r})dL_{r[O]}.
\end{align}
Note that $\int_{G_{r,n-r}}V_{r}(\Phi^{'}_{r})dL_{r[O]}=I_{n-r}(\Phi)$ and $W^{(n)}_{n-r}(\Phi)=
\frac{rO_{r-2}\cdots O_{0}}{nO_{n-2}\cdots O_{n-r}}I_{n-r}(\Phi)$, then we have
\begin{equation}\label{d2}
\int_{G_{r,n-r}}V_{r}(\Phi^{'}_{r})dL_{r[O]}=\frac{nO_{n-2}\cdots O_{n-r}}{rO_{r-2}\cdots O_{0}}W^{(n)}_{n-r}(\Phi)=\frac{O_{n-2}\cdots O_{n-r}}{rO_{r-2}\cdots O_{0}}M^{(n)}_{n-r-1}(\partial\Phi).
\end{equation}
Inserting \eqref{d2} to \eqref{d1} and using \eqref{c4}, we get
\begin{align}\label{d3}
\int_{G_{r,n-r}}M^{(n)}_{l}&\big(\partial(\Phi^{'}_{r})^{(n)}_{\rho}\big)dL_{r[O]}=\notag\\
&(-1)^{r}\tbinom{n-1}{n-r-1}^{-1}\frac{O_{n-r-1}O_{n-2}\cdots O_{n-r}}{rO_{r-2}\cdots O_{0}}M^{(n)}_{n-r-1}(\partial\Phi)\notag\\
&+\sum^{r-1}_{i=0}(-1)^{i}
 \frac{\tbinom{r}{i}\tbinom{r-1}{r-i-1}}{\tbinom{n-1}{n-i-1}}\tfrac{O_{n-i-1}O_{n-2}\cdots O_{n-r}}{O_{r-i-1}O_{r-2}\cdots O_{0}}h^{r-i}M^{(n)}_{n-i-1}(\partial\Phi)\notag\\
 &+\sum^{r}_{j=1}\sum^{r-j}_{i=0}(-1)^{i}\frac{\tbinom{r}{j}\tbinom{r-j}{i}\tbinom{r-1}{r-i-1}}{\tbinom{n-1}{n-i-1}}
 \tfrac{O_{n-i-1}O_{n-2}\cdots O_{n-r}}{O_{r-i-1}O_{r-2}\cdots O_{0}}\rho^{j}h^{r-j-i}M^{(n)}_{n-i-1}(\partial\Phi).
\end{align}

\textbf{Case 3.}When $q<n-r-1$, by \eqref{case3},\eqref{c4} and \eqref{d2}, we have
\begin{align*}
&\int_{G_{r,n-r}}M^{(n)}_{l}\big(\partial(\Phi^{'}_{r})^{(n)}_{\rho}\big)dL_{r[O]}\notag=\\
  &\sum^{n-r-l-1}_{j=0}(-1)^{r}
  \tbinom{n-l-1}{j}
  \tbinom{n-l-j-1}{r}
  \tbinom{n-1}{n-r-1}^{-1}O_{n-r-1}\rho^{j}h^{n-l-j-r-1}\int_{G_{r,n-r}}V_{r}(\Phi^{'}_{r})dL_{r[O]}\notag\\
  &+\sum^{n-r-l}_{j=0}\sum^{r-1}_{i=0}(-1)^{i}
  \frac{\tbinom{n-l-1}{j}\tbinom{n-l-j-1}{i}\tbinom{r-1}{r-i-1}}{\tbinom{n-1}{n-i-1}}\tfrac{O_{n-i-1}}{O_{r-i-1}}\rho^{j}h^{n-l-j-i-1}
 \int_{G_{r,n-r}}M^{(r)}_{r-i-1}(\partial\Phi^{'}_{r})dL_{r[O]}\notag\\
 &+\sum^{n-l-1}_{j=n-r-l+1}\sum^{n-l-j-1}_{i=0}(-1)^{i}
  \frac{\tbinom{n-l-1}{j}\tbinom{n-l-j-1}{i}\tbinom{r-1}{r-i-1}}{\tbinom{n-1}{n-i-1}}\tfrac{O_{n-i-1}}{O_{r-i-1}}\rho^{j}h^{n-l-j-i-1}
 \int_{G_{r,n-r}}M^{(r)}_{r-i-1}(\partial\Phi^{'}_{r})dL_{r[O]}\notag\\
 &=\sum^{n-r-l-1}_{j=0}(-1)^{r}
 \tfrac{\tbinom{n-l-1}{j}\tbinom{n-l-j-1}{r}}{\tbinom{n-1}{n-r-1}}\tfrac{O_{n-r-1}O_{n-2}\cdots O_{n-r}}{rO_{r-2}\cdots O_{0}}\rho^{j}h^{n-l-j-r-1}M^{(n)}_{n-r-1}(\partial\Phi)\notag\\
  &+\sum^{n-r-l}_{j=0}\sum^{r-1}_{i=0}(-1)^{i}
  \frac{\tbinom{n-l-1}{j}\tbinom{n-l-j-1}{i}\tbinom{r-1}{r-i-1}}{\tbinom{n-1}{n-i-1}}\tfrac{O_{n-i-1}O_{n-2}\cdots O_{n-r}}{O_{r-i-1}O_{r-2}\cdots O_{0}}\rho^{j}h^{n-l-j-i-1}
 M^{(n)}_{n-i-1}(\partial\Phi)\notag\\
 &+\sum^{n-l-1}_{j=n-r-l+1}\sum^{n-l-j-1}_{i=0}(-1)^{i}
  \frac{\tbinom{n-l-1}{j}\tbinom{n-l-j-1}{i}\tbinom{r-1}{r-i-1}}{\tbinom{n-1}{n-i-1}}\tfrac{O_{n-i-1}O_{n-2}\cdots O_{n-r}}{O_{r-i-1}O_{r-2}\cdots O_{0}}\rho^{j}h^{n-l-j-i-1}
 M^{(n)}_{n-i-1}(\partial\Phi).
\end{align*}
 This completes the proof Theorem \eqref{them2}.
\end{proof}

\textbf{Acknowledgements}
This work was partially supported by Science and Technology Commission
of Shanghai Municipality (No. 22DZ2229014) and the National Natural Science Foundation of
China (Grant No. 12271163). The research is supported by Shanghai Key Laboratory of PMMP.

\textbf{Data availability}
Data sharing not applicable to this article as no datasets were generated or analysed during the current study.

\textbf{Conflict of interest}
The author has no conflicts of interest to declare that are relevant to the content of this article.

\end{document}